\documentclass[12pt,reqno,a4paper]{amsart}
\usepackage{latexsym,bm}
\usepackage{amsmath,amssymb,cases}
\usepackage{amsthm}
\usepackage{xcolor}
\usepackage{enumerate}
\usepackage{caption}
\usepackage{graphicx,subfig}
\usepackage[numbers,sort&compress]{natbib}
\usepackage{accents}

\usepackage{geometry}
\usepackage{float}

\usepackage{ulem}
\usepackage[thicklines]{cancel}
\usepackage{relsize}
\usepackage{exscale}

\newtheorem{The}{Theorem}[section]
\newtheorem{Lem}[The]{Lemma}

\newtheorem{Def}[The]{Definition}

\theoremstyle{definition}

\graphicspath{{fig/}}

\numberwithin{equation}{section}
\numberwithin{figure}{section}

\begin{document}
\captionsetup[figure]{labelfont={default},labelformat={default},labelsep=period,name={Fig.}}

\title[Zero--relaxation and vanishing--damping limits of pressureless Euler system]
{Zero--relaxation and vanishing--damping limits of pressureless Euler system
}

\author{Guirong Tang}
\address{Guirong Tang: School of Mathematical Sciences, Capital Normal University, Beijing 100048, China}
\email{\tt tangguirong@amss.ac.cn}

\begin{abstract}
We are concerned with the one-dimensional pressureless Euler system with relaxation in the Radon measure space.
As the relaxation time tends to zero, the entropy solution converges to a static solution with the density converging to its initial value.
As the relaxation time tends to infinity, which means the damping vanishes,
the entropy solution of damped pressureless Euler system converges to that of pressureless Euler system.
\end{abstract}

\date{\today}
\maketitle

\tableofcontents
 
\thispagestyle{empty}

\section{Introduction}
We are concerned with the following one-dimensional pressureless Euler system with relaxation/damping term:
\begin{equation}\label{1.1}
\begin{cases}
\rho_t +(\rho u) _x=0,\\[1mm]
(\rho u)_t+(\rho u^2)_x=-\frac{\rho u}{\tau},
\end{cases}
\end{equation}
for $ (x,t) \in \mathbb{R}^2_+:=({-}\infty,\infty)\times(0,\infty)$.   
Here $\rho\geq 0$ is the density of mass, $u$ denotes the velocity,
$\tau>0$ is the relaxation time which means the average time required for the state being relaxed to the equilibrium.   

\vspace{1pt}
Since system $\eqref{1.1}$ is linearly degenerate hyperbolic, 
then the solutions develop singularities and form $\delta$-shocks generically in a finite time for a large class of intial data.
Thus, it is natural and necessary to understand the solutions in the sense of Radon measures.
We consider the general Cauchy initial data as follows:
\begin{equation}\label{1.2}
(\rho,u)|_{t=0}=(\rho_0,u_0) \in (\mathcal{M}_{loc}(\mathbb{R}),L^{\infty}_{\rho_0}(\mathbb{R})),
\end{equation}
where $\mathcal{M}_{loc}(\mathbb{R})$ is the space of locally finite Radon measures
and $L^{\infty}_{\rho_0}(\mathbb{R})$ represents the space of bounded measurable functions 
with respect to the measure $\rho_0$.

The existence and uniqueness of entropy solution to \eqref{1.1} have been studied by \cite{HHW2014} and \cite{Jin2016} with $\tau=1$.
In \cite{HHW2014}, Ha-Huang-Wang consider the initial density being a Lebesgue-measurable function and get the formula of entropy solution.
Jin \cite{Jin2016} extends their result to the case that the initial density is locally finite Radon measure.
When the damping term is multiplied by a coefficient $\frac{1}{\tau}$, their results remain valid, requiring only a slight modification to the formulation, while the proof procedure remains essentially identical.

\smallskip
As the damping term is vanishing, {\it i.e.} $\tau\rightarrow\infty$, 
if denoting the limits of $(\rho,u)$ as $(\bar{\rho},\bar{u})$,
the system \eqref{1.1} formally transforms into pressureless Euler system:
\begin{align}\label{1.3}
\begin{cases}
\bar{\rho}_t+(\bar{\rho} \bar{u})_x=0,\\[1mm]
(\bar{\rho}\bar{u})_t+(\bar{\rho}\bar{u}^2)_x=0.
\end{cases}
\end{align}
It is well-known that the existence and uniqueness of entropy solution to \eqref{1.3} have been obtained by Huang-Wang \cite{huang2001well}.
In this paper, we rigorously prove that the entropy solution to \eqref{1.1} converges to the entropy solution to \eqref{1.3} as the damping term vanishes.

\smallskip
As the relaxation time tends to zero, {\it i.e.} $\tau\rightarrow0{+}$, by slow time scaling \cite{MR2000}, {\it i.e.} $t=\tau t'$,
if we introduce 
\begin{align}\label{1.3b}
\rho^{\tau}(x,t)=\rho(x,\frac{t}{\tau}),\quad u^{\tau}(x,t)=\frac{1}{\tau}u(x,\frac{t}{\tau}),
\end{align}
then the degenerate hyperbolic system \eqref{1.1}--\eqref{1.2} transforms into
\begin{equation}\label{1.4}
\begin{cases}
\rho^{\tau}_t +(\rho^{\tau} u^{\tau})_x=0,\\[1mm]
\tau^2\big((\rho^{\tau} u^{\tau})_t+(\rho^{\tau} (u^{\tau})^2)_x \big)=-\rho^{\tau} u^{\tau}.
\end{cases}
\end{equation}
with initial data
\begin{align}
\rho^{\tau}_0(x):=\rho^{\tau}(x,0)=\rho_0(x), \qquad
u^{\tau}_0(x):=u^{\tau}(x,0)=\frac{1}{\tau} u_0(x).
\end{align}
Formally, as $\tau \rightarrow 0$, if denoting the limits of $(\rho^{\tau},  u^{\tau})$ as $(\hat{\rho},\hat{u})$, 
one gets that the formal limiting system of \eqref{1.4} is
\begin{align}\label{1.5}
\begin{cases}
\hat{\rho}_t+(\hat{\rho} \hat{u})_x=0,\\[1mm]
\hat{\rho}\hat{u}=0,
\end{cases}
\end{align}
which formally implies that $\hat{\rho}(x,t)=\rho_0(x)$ and $\hat{u}(x,t)\equiv 0$ on the support of $\hat{\rho}$.
In this paper, we rigorously show that the relaxation limit of the entropy solution to \eqref{1.1} is $(\rho_0,0)$ by slow time scaling.

There are many physical models with relaxation time and the relaxation limit is well-known in asymptotic analysis and singular perturbation theories \cite{GBW}.
The limit can be from first-order hyperbolic systems to hyperbolic system (see \cite{Liu1987, Na1997} and references therein), 
or to parabolic system by slow time scaling \cite{MR2000}.

The main results concentrate on the damped Euler system with pressure.
Junca--Rascle \cite{JR2002} consider large BV solution of isothermal Euler equations in one dimension for the initial data being $L^1$ perturbations of Riemann data and prove that the density converges to the solution of the heat equation.
The result is extended to multidimensional case by Coulombel--Goudon \cite{CG2007}.
For polytropic gas in one dimension, Marcati--Milani \cite{MM1990}
prove the limit density satisfies Darcy's law, so that it is a weak solution of the porous media.
The convergence can also be considered in Besov space.
Xu \cite{X2010} consider the relaxation limit of  multidimensional isentropic Euler equations in the framework of Besov spaces with relatively lower regularity.
Xu--Wang \cite{XW2013} consider the multidimensional Euler equations in the Chemin--Lerner's space with initial data close to the constant equilibrium state.
Xu--Kawashima \cite{XK2014} obtain the relaxation limit of Euler equations towards the porous medium equation by Aubin--Lions compactness argument.
For non-isentropic Euler equations in $\mathbb{R}^3$, Wu \cite{W2016} prove the existence of initial layer for ill-prepared data and get the strong convergence rates.
Lattanzio--Tzavaras \cite{LT2013} establish the relaxation limit from damped Euler equations to porous media equation away from vacuum by relative entropy method. 
Furthermore, there are many studies about the relaxation limit of Euler--Poisson system, such as \cite{MN1995, ABR2000, HZ2000, Wang2001, LM1999, LPZ2021}.

There are few results about the relaxation limit of linearly degenerated system.
For damped pressureless Euler system, under the sole influence of damping, the particle's initial kinetic energy dissipates over time, resulting in asymptotic decay of its velocity to zero at equilibrium.
In the limit of vanishing damping, the particle's dynamics reduce to undamped motion with constant initial velocity, thereby recovering the solution of the pressureless Euler system.
Now we state the main theorem of this paper.
\begin{The}
Let $(\rho(x,t),u(x,t))$ and $(\bar{\rho}(x,t),\bar{u}(x,t))$ be the unique entropy solution of \eqref{1.1} and \eqref{1.3} respectively with the initial data $(\rho_0,u_0)$ in \eqref{1.2},
let $(\rho^{\tau}(x,t),u^{\tau}(x,t))$ be defined by \eqref{1.3b}.
Then,
\begin{align}
&\rho^{\tau}(x,t)\rightharpoonup \rho_0 \qquad\; {\rm and}\quad u^{\tau}(x,t)\rightarrow 0, \qquad \quad \;\,\,{\rm as}\ \tau\rightarrow 0;
\\[1mm]
&\rho(x,t)\rightharpoonup \bar{\rho}(x,t) \quad {\rm and}\quad u(x,t)\rightarrow \bar{u}(x,t), \qquad {\rm as}\ \tau\rightarrow \infty.
\end{align}
\end{The}
\smallskip

\vspace{2pt}
This paper is organized as follows: 
In \S 2, we introduce the formula of entropy solutions and some relevant results about pressureless Euler system with or without damping.
In \S 3, we rigorously prove the zero--relaxation and vanishing--damping limits.

\section{Formula of entropy solutions}

In this section, we introduce the formula of entropy solutions and some relevant results about pressureless Euler system with/without damping in Jin \cite{Jin2016} and Huang-Wang \cite{huang2001well}.
In Jin \cite{Jin2016}, the damping term is ${-}\rho u$ without the relaxation time $\tau$. 
We will consider this form ${-}\frac{\rho u}{\tau}$ and restate the contents without proof, since the arguments are almost the same.

\subsection{Pressureless Euler system.}
In Huang-Wang \cite{huang2001well}, by introducing a function 
\begin{align}
\bar{m}(x,t)=\oint_{(0,0)}^{(x,t)} \bar{\rho}\,{\rm d}x{-}\bar{\rho} \bar{u}\,{\rm d}t,
\end{align}
the pressureless Euler system transforms into
\begin{align}\label{2.1}
\begin{cases}
\bar{m}_t+\bar{u} \bar{m}_x=0,\\[1mm]
(\bar{m}_x \bar{u})_t+(\bar{m}_x \bar{u}^2)_x=0.
\end{cases}
\end{align}
The definitions of weak solution and entropy solution are as follows.

\begin{Def}[Weak solution I]\label{Def1.1}
Let $\bar{m}(x,t)$ be of bounded variation locally in $x$ and $\bar{u}(x,t)$ be bounded and measurable to $\bar{m}_x$.
Assume that the measures $\bar{m}_x$ and $\bar{u} \bar{m}_x$ are weakly continuous in $t$.
$(\bar{\rho},\bar{u})$ is called a weak solution of \eqref{1.3} or $(\bar{m},\bar{u})$ is called a weak solution of \eqref{2.1}, if
\begin{align}
\begin{cases}
\iint \varphi_t \bar{m}\,{\rm d}x {\rm d}t-\iint\varphi \bar{u}\,{\rm d}\bar{m} {\rm d}t=0,\\[1mm]
\iint \psi_t \bar{u}+\psi_x \bar{u}^2\,{\rm d}\bar{m} {\rm d}t=0,
\end{cases}
\end{align}
holds for all $\varphi, \psi\in C_0^{\infty}(\mathbb{R}^2_{+})$.
Here $\iint \cdots \,{\rm d}\bar{m}{\rm d}t$ denotes Lebesgue-Stieltjes integral.

The initial value is understood in the following sense:
as $t\rightarrow 0{+}$, the measures $\bar{\rho}$ and $\bar{\rho} \bar{u}$ weakly converge to $\rho_0$, $\rho_0 u_0$ respectively.
\end{Def}

\begin{Def}[Entropy solution I]\label{Def1.2}
Let $(\bar{\rho},\bar{u})$ be a weak solution of \eqref{1.3}.
$(\bar{\rho},\bar{u})$ is called an entropy solution of \eqref{1.3} if 
\begin{align}
\frac{\bar{u}(x_2,t)-\bar{u}(x_1,t)}{x_2-x_1}\leq \frac{1}{t},
\end{align}
holds for any $x_1<x_2$, almost everywhere $t>0$ and the measure $\bar{\rho} \bar{u}^2$ weakly converges to $\rho_0 u_0^2$ as $t\rightarrow 0$.
\end{Def}

Now we state the construction of the formula for the unique entropy solution.
Firstly, introduce the generalized potential 
\begin{align}\label{2.4}
\bar{F}(y;x,t)=\int_{0{+}0}^{y{-}0} \eta+t u_0(\eta)-x\,{\rm d}m_0(\eta),
\end{align}
where $m_0(x)=\rho_0([0,x[)$.
Define the minimum value and the minimum points for fixed $(x,t)$ as follows:
\begin{align}
&\bar{\nu}(x,t)=\min_y \bar{F}(y;x,t),\label{2.5}\\[1mm]
&\bar{S}(x,t)=\{y; \exists y_n\rightarrow y,\ {\rm s.t.}\ \bar{F}(y_n;x,t)\rightarrow \bar{\nu}(x,t)\}. \label{2.6}
\end{align}
If $y_0\in \bar{S}(x,t)$ with $[m_0(y_0)]=m_0(y_0{+}0)-m_0(y_0{-}0)>0$, then
\begin{align}
\bar{\nu}(x,t)=
\begin{cases}
\bar{F}(y_0;x,t) \qquad &{\rm if}\ x\leq y_0+t u_0(y_0);\\[1mm]
\bar{F}(y_0{+}0;x,t) \qquad &{\rm if}\ x> y_0+t u_0(y_0).
\end{cases}
\end{align}
Let 
\begin{align}\label{2.8}
y_m(x,t)=\inf_{y\in \bar{S}(x,t)}y, \qquad 
y^m(x,t)=\sup_{y\in \bar{S}(x,t)}y,
\end{align}
and 
\begin{align*}
y_M(x,t)=
\begin{cases}
y_{\#}, \qquad {\rm if}\ (x,t)\in A(x,t)\ {\rm and}\ y_{\#}\in {\rm spt}\{\rho_0\}, \\[1mm]
y_m, \qquad {\rm otherwise,}
\end{cases}\\[2mm]
y^M(x,t)=
\begin{cases}
y^{\#}, \qquad {\rm if}\ (x,t)\in B(x,t)\ {\rm and}\ y^{\#}\in {\rm spt}\{\rho_0\}, \\[1mm]
y^m, \qquad {\rm otherwise,}
\end{cases}
\end{align*}
where 
\begin{align*}
A(x,t)=\{(x,t); \exists y_{\#}>y_m,\ {\rm s.t.}\ \rho_0=0 \ {\rm on}\ (y_m,y_{\#})\ {\rm and}\ \forall y \in (y_m,y_{\#}),\ y\in \bar{S}(x,t)\},\\[1mm]
B(x,t)=\{(x,t); \exists y^{\#}<y^m,\ {\rm s.t.}\ \rho_0=0 \ {\rm on}\ (y^{\#},y^m)\ {\rm and}\ \forall y \in (y^{\#},y^m),\ y\in \bar{S}(x,t)\}.
\end{align*}
We take two special minimum points in $\bar{S}(x,t)\cap {\rm spt}\{\rho_0\}$,
\begin{align}\label{2.9}
\bar{y}_*(x,t)=\min\{y_M(x,t), y^M(x,t)\}, \qquad \bar{y}^*(x,t)=y^M(x,t).
\end{align}
Then $\bar{y}_*(x,t)$ and $\bar{y}^*(x,t)$ are increasingly monototic in $x$.
In particular, for fixed $t>0$ and $x_1<x_2$, 
\begin{align}
\bar{y}_*(x_1,t)\leq \bar{y}^*(x_2,t),
\end{align}
which implies that, for fixed $t$ and almost all $x\in \mathbb{R}$,
\begin{align}\label{2.10}
\bar{y}_*(x,t)=\bar{y}^*(x,t).
\end{align}
The mass and momentum of pressureless Euler system can be defined as follows:
\begin{align}
&\bar{m}(x,t)=
\begin{cases}
\int_{0{+}0}^{\bar{y}_*(x,t){-}0}\,{\rm d}m_0(\eta), \quad{\rm if}\ \bar{\nu}(x,t)=\bar{F}(\bar{y}_*(x,t);x,t),\\[1mm]
\int_{0{+}0}^{\bar{y}_*(x,t){+}0}\,{\rm d}m_0(\eta), \quad{\rm if}\ \bar{\nu}(x,t)=\bar{F}(\bar{y}_*(x,t){+}0;x,t);
\end{cases}\\[1mm]
&\bar{q}(x,t)=
\begin{cases}
\int_{0{+}0}^{\bar{y}_*(x,t){-}0}u_0(\eta)\,{\rm d}m_0(\eta), \quad{\rm if}\ \bar{\nu}(x,t)=\bar{F}(\bar{y}_*(x,t);x,t),\\[1mm]
\int_{0{+}0}^{\bar{y}_*(x,t){+}0}u_0(\eta)\,{\rm d}m_0(\eta), \quad{\rm if}\ \bar{\nu}(x,t)=\bar{F}(\bar{y}_*(x,t){+}0;x,t).
\end{cases}
\end{align}

For almost everywhere $t>0$, the velocity $\bar{u}(x,t)$ satisfies
\begin{align}\label{2.11}
\bar{q}_x(x,t)=\bar{u}(x,t) \bar{m}_x(x,t),
\end{align}
in the sense of Radon-Nikodym derivative.

\subsection{Damped pressureless Euler system.}
In Jin \cite{Jin2016}, by introducing a function
\begin{align}
m(x,t)=\oint_{(0,0)}^{(x,t)} \rho\,{\rm d}x{-}\rho u\,{\rm d}t,
\end{align}
the damped pressureless Euler system \eqref{1.1} transforms into
\begin{align}
\begin{cases}
m_t+u m_x=0,\\[1mm]
(m_x u)_t+(m_x u^2)_x=-\frac{m_x u}{\tau}.
\end{cases}
\end{align}
The definitions of weak solution and entropy solution are similar with Definition \ref{Def1.1} and Definition \ref{Def1.2}.

\begin{Def}[Weak solution II]\label{Def2.1}
Let $m(x,t)$ be of bounded variation locally in $x$ and $u(x,t)$ be bounded and measurable to $m_x$.
Assume that the measures $m_x$ and $u m_x$ are weakly continuous in $t$.
$(\rho,u)$ is called a weak solution of \eqref{1.1} or $(m,u)$ is called a weak solution of \eqref{2.1}, if
\begin{align}
\begin{cases}
\iint \varphi_t m\,{\rm d}x {\rm d}t-\iint\varphi u\,{\rm d}m {\rm d}t=0,\\[1mm]
\iint \psi_t u+\psi_x u^2- \frac{1}{\tau}\psi u\,{\rm d}m {\rm d}t=0,
\end{cases}
\end{align}
holds for all $\varphi, \psi\in C_0^{\infty}(\mathbb{R}^2_{+})$.
Here $\iint \cdots \,{\rm d}m{\rm d}t$ denotes Lebesgue-Stieltjes integral.

The initial value is understood in the following sense:
as $t\rightarrow 0{+}$, the measures $\rho$ and $\rho u$ weakly converge to $\rho_0$, $\rho_0 u_0$ respectively.
\end{Def}

\begin{Def}[Entropy solution II]\label{Def2.2}
Let $(\rho,u)$ be a weak solution of \eqref{1.1}.
$(\rho,u)$ is called an entropy solution of \eqref{1.1} if 
\begin{align}\label{2.18}
\frac{u(x_2,t)-u(x_1,t)}{x_2-x_1}\leq \frac{e^{-\frac{t}{\tau}}}{\tau(1-e^{-\frac{t}{\tau}})},
\end{align}
holds for any $x_1<x_2$, almost everywhere $t>0$ and the measure $\rho u^2$ weakly converges to $\rho_0 u_0^2$ as $t\rightarrow 0$.
\end{Def}

If we consider the smooth solution, the momentum equation is
\begin{align*}
u_t+u u_x=-\frac{1}{\tau} u,
\end{align*}
The trajectory $X(\eta,t)$ of a particle originated from $\eta$ with initial velocity $u_0(\eta)$ satisfies the following ODE:
\begin{align*}
\begin{cases}
\frac{{\rm d}X(\eta,t)}{{\rm d}t}=u(X(\eta,t),t),\\[1mm]
\frac{{\rm d}u(X(\eta,t),t)}{{\rm d}t}=-\frac{1}{\tau}u(X(\eta,t),t),\\[1mm]
X(\eta,0)=\eta, \quad u(\eta,0)=u_0(\eta),
\end{cases}
\end{align*}
which implies that
\begin{align*}
X(\eta,t)=\eta+u_0(\eta)(\tau-\tau e^{-\frac{t}{\tau}}).
\end{align*}
Then the generalized potential in Jin \cite{Jin2016} is adjusted to this form:
\begin{align}\label{2.19}
F(y;x,t)=\int_{0{+}0}^{y{-}0}\eta+u_0(\eta)(\tau-\tau e^{-\frac{t}{\tau}})-x\,{\rm d}m_0(\eta).
\end{align}
Similar to \eqref{2.5} and \eqref{2.6},
\begin{align}
&\nu(x,t)=\min_{y}F(y;x,t),\\[1mm]
&S(x,t)=\{y; \exists y_n\rightarrow y,\ {\rm s.t.}\ F(y_n;x,t)\rightarrow \nu(x,t)\}.
\end{align}
If $y_0\in S(x,t)$ with $[m_0(y_0)]=m_0(y_0{+}0)-m_0(y_0{-}0)>0$, then
\begin{align}
\nu(x,t)=
\begin{cases}
F(y_0;x,t) \qquad &{\rm if}\ x\leq y_0+ u_0(y_0)(\tau-\tau e^{-\frac{t}{\tau}});\\[1mm]
F(y_0{+}0;x,t) \qquad &{\rm if}\ x> y_0+u_0(y_0)(\tau-\tau e^{-\frac{t}{\tau}}).
\end{cases}
\end{align}
By the same procedures as \eqref{2.8}--\eqref{2.9}, we can choose two special minimum points in $S(x,t)\cap {\rm spt}\{\rho_0\}$:
\begin{align}\label{2.22}
y_*(x,t)\quad {\rm and}\quad y^*(x,t).
\end{align}
Then $y_*(x,t)$ and $y^*(x,t)$ are increasingly monototic in $x$.
In particular, for fixed $t>0$ and $x_1<x_2$, 
\begin{align}
y_*(x_1,t)\leq y^*(x_2,t),
\end{align}
which implies that, for fixed $t$ and almost all $x\in \mathbb{R}$,
\begin{align}\label{2.20}
y_*(x,t)=y^*(x,t).
\end{align}
For $(x,t)$ satisfying \eqref{2.20}, denote $U=||u_0||_{L^{\infty}}$, then 
\begin{align}\label{2.21}
x-U(\tau-\tau e^{-\frac{t}{\tau}})\leq y_*(x,t) \leq x+U(\tau-\tau e^{-\frac{t}{\tau}}).
\end{align}

For each point $(x_0,t_0)$, the left and right backward characteristics $L_1, L_2$ are defined as follows:
\begin{align*}
L_1: \frac{x-y_*(x_0,t_0)}{\tau-\tau e^{-\frac{t}{\tau}}}=\frac{x_0-y_*(x_0,t_0)}{\tau-\tau e^{-\frac{t_0}{\tau}}},\\[1mm]
L_2: \frac{x-y^*(x_0,t_0)}{\tau-\tau e^{-\frac{t}{\tau}}}=\frac{x_0-y^*(x_0,t_0)}{\tau-\tau e^{-\frac{t_0}{\tau}}}.
\end{align*}
$L_1, L_2$ and $x$-axis form a area and we denote it as $\Delta(x_0,t_0)$.

The mass and the momentum of damped pressureless Euler system are defined as follows:
\begin{align}
&m(x,t)=
\begin{cases}
\int_{0{+}0}^{y_*(x,t){-}0}\,{\rm d}m_0(\eta), \quad{\rm if}\ \nu(x,t)=F(y_*(x,t);x,t),\\[1mm]
\int_{0{+}0}^{y_*(x,t){+}0}\,{\rm d}m_0(\eta), \quad{\rm if}\ \nu(x,t)=F(y_*(x,t){+}0;x,t);
\end{cases}\\[1mm]
&q(x,t)=
\begin{cases}
\int_{0{+}0}^{y_*(x,t){-}0}u_0(\eta)e^{-\frac{t}{\tau}}\,{\rm d}m_0(\eta), \quad{\rm if}\ \nu(x,t)=F(y_*(x,t);x,t),\\[1mm]
\int_{0{+}0}^{y_*(x,t){+}0}u_0(\eta)e^{-\frac{t}{\tau}}\,{\rm d}m_0(\eta), \quad{\rm if}\ \nu(x,t)=F(y_*(x,t){+}0;x,t).
\end{cases}\label{2.21b}
\end{align}

Now we give the formula of velocity $u(x,t)$ without proof when the damping is $-\frac{\rho u}{\tau}$,
which is presented by Lemma 2.6 in \cite{Jin2016}.

\begin{Lem}
Each point $(x_0,t_0)$ at $t_0>0$ uniquely determines a Lipschitz continuous curve $L$: $x=x(t)$ with $x_0=x(t_0)$.
In adddition, for any $t\in\{\tau: \tau\geq t_0\}$,
\begin{align*}
x'(t)=\lim_{t'',t'\rightarrow t{+}0}\frac{x(t'')-x(t')}{t''-t'}=
\begin{cases}
\frac{x-y_*(x,t)}{\tau-\tau e^{-\frac{t}{\tau}}}e^{-\frac{t}{\tau}}, &(x,t)\in V_3\cap V_4,\\[1mm]
\frac{x-y_*(x{+}0,t)}{\tau-\tau e^{-\frac{t}{\tau}}}e^{-\frac{t}{\tau}}, &(x,t)\in V_5,\\[1mm]
\lim_{x_2, x_1\rightarrow x{\pm}0}
\frac{\int_{\tilde{y}_*(x_1,t){-}0}^{\tilde{y}^*(x_2,t){+}0}u_0(\eta)\,{\rm d}m_0(\eta)}
{\int_{\tilde{y}_*(x_1,t){-}0}^{\tilde{y}^*(x_2,t){+}0}\,{\rm d}m_0(\eta)}e^{-\frac{t}{\tau}}, &(x,t)\in V_1\cap V_2,
\end{cases}
\end{align*}
where
\begin{align*}
&V_1=\{(x,t)\in L: y_*(x,t)<y^*(x,t)\},\\[1mm]
&V_2=\{(x,t)\in L: y_*(x,t)=y^*(x,t), \forall t'>t, y_*(x(t'),t')<y^*(x(t'),t')\\[1mm] 
&\qquad\qquad\qquad\qquad {\rm and}\ y_*(x{-}0,t)<y_*(x{+}0,t)\}\\[1mm]
&V_3=\{(x,t)\in L:y_*(x,t)=y^*(x,t), \forall t'>t, y_*(x(t'),t')<y^*(x(t'),t')\\[1mm] 
&\qquad\qquad\qquad\qquad {\rm and}\ y_*(x{-}0,t)=y_*(x{+}0,t)\}\\[1mm]
&V_4=\{(x,t)\in L: y_*(x,t)=y^*(x,t), \exists t_a>t\ {\rm s.t.}\ y_*(x(t_a),t_a)= y^*(x(t_a),t_a)\\[1mm]
&\qquad\qquad\qquad\qquad {\rm and}\ (x,t)\in \Delta(x(t_a),t_a)\},\\[1mm]
&V_4=\{(x,t)\in L: y_*(x,t)=y^*(x,t), \exists t_a>t\ {\rm s.t.}\ y_*(x(t_a),t_a)= y^*(x(t_a),t_a)\\[1mm]
&\qquad\qquad\qquad\qquad {\rm and}\ (x,t)\notin \Delta(x(t_a),t_a)\},\\[1mm]
&\tilde{y}_*(x_1,t)=y_*(x_1,t){+}(x-x_1) H\big(x_1-(y_*(x{-}0,t)+u_0(y_*(x{-}0,t)(\tau-\tau e^{-\frac{t}{\tau}})\big)\\[1mm]
&\tilde{y}^*(x_2,t)=y^*(x_2,t){+}(x-x_2)H\big(y_*(x{+}0,t)+u_0(y_*(x{+}0,t)(\tau-\tau e^{-\frac{t}{\tau}})-x_2\big),\\[1mm]
&H(x)=
\begin{cases}
0,\qquad x\leq 0,\\[1mm]
1,\qquad x>0.
\end{cases}
\end{align*}
\end{Lem}

By making adjustment in the vacuum region caused by initial vacuum instead of rarefaction, 
the velocity $u(x,t)$ is defined as:
\begin{align}
u(x,t)=
\begin{cases}
-U e^{-\frac{t}{\tau}}, &{\rm if}\ x-y_*(x,t)<-U (\tau-\tau e^{-\frac{t}{\tau}}),\\[1mm]
x'(t), &{\rm if}\ |x-y_*(x,t)|\leq U (\tau-\tau e^{-\frac{t}{\tau}}),\\[1mm]
U e^{-\frac{t}{\tau}}, &{\rm if}\ x-y_*(x,t)>U (\tau-\tau e^{-\frac{t}{\tau}}),
\end{cases}
\end{align}
where $U=||u_0||_{L^{\infty}}$, $x(t)$ is the Lipschitz continuous curve in above lemma.

For $(x,t)$ satisfies $|x-y_*(x,t)|\leq U (\tau-\tau e^{-\frac{t}{\tau}})$, we have
\begin{align}
&u(x{-}0,t)=\frac{x-y_*(x{-}0,t)}{\tau-\tau e^{-\frac{t}{\tau}}}e^{-\frac{t}{\tau}},\quad
u(x{+}0,t)=\frac{x-y_*(x{+}0,t)}{\tau-\tau e^{-\frac{t}{\tau}}}e^{-\frac{t}{\tau}},\\[1mm]
&u(x{+}0,t)\leq u(x,t)\leq u(x{-}0,t),
\end{align}
which implies the Oleinik type entropy condition \eqref{2.18}.

For almost everywhere $t>0$, the velocity $u(x,t)$ satisfies
\begin{align}\label{2.23}
q_x(x,t)=u(x,t) m_x(x,t),
\end{align}
in the sense of Radon-Nikodym derivative.

\section{Zero--relaxation and Vanishing--damping limits}

This section is devoted to the proof of Theorem 1.1.
The strategy is as follows: Firstly, we prove the convergence of generalized potential, which yields the minimum points have convergent subsequence.
Secondly, by the formula of mass, we get the weak convergence of density.
The limit of the velocity is follows from the convergence of momentum and the Radon-Nikodym derivative.

\subsection{Zero--relaxation limit.}
By slow time scaling, {\it i.e.} $t=\tau t'$, the generalized potential $F(y;x,t)$ transforms into
\begin{align}
F^{\tau}(y;x,t):=F(y;x,\frac{t}{\tau})=\int_{0{+}0}^{y{-}0}\eta+u_0(\eta)(\tau-\tau e^{-\frac{t}{\tau^2}})-x\,{\rm d}m_0(\eta),
\end{align}
and the minimum point in \eqref{2.22} transforms into
\begin{align}
y^{\tau}_*(x,t):=y_*(x,\frac{t}{\tau}).
\end{align}
For $(x,t)$ satisfying $y_*(x,t)=y^*(x,t)$, it follows from \eqref{2.21} that $\{y^\tau(x,t)\}$ is a bounded sequence with respect to $\tau$.
Thus there exists a convergent subsequence (still denoted as $y^\tau(x,t)$) with
\begin{align}
\hat{y}(x,t):=\lim_{\tau\rightarrow 0}y^\tau(x,t)\in {\rm spt}\{\rho_0\}.
\end{align}
As the relaxation time $\tau\rightarrow 0$, it is direct to see that 
\begin{align}\label{3.1}
\hat{F}(y;x,t):=\lim_{\tau\rightarrow 0}F^{\tau}(y;x,t)=\int_{0{+}0}^{y{-}0} \eta-x\,{\rm d}m_0(\eta).
\end{align}
Since $F^{\tau}(y_*^{\tau}(x,t);x,t)\leq F^{\tau}(y;x,t)$ for any $y\in \mathbb{R}$, by the left continuity of $F^{\tau}(y;x,t)$ with respect to $y$, at least one of the following holds:  
\begin{align}\label{3.1a}
\hat{F}(\hat{y}(x,t);x,t)\leq \hat{F}(y;x,t), \quad {\rm or}\quad 
\hat{F}(\hat{y}(x,t){+}0;x,t)\leq \hat{F}(y;x,t).
\end{align}
From \eqref{3.1}, obviously for $(x,t)$ satisfying $y_*(x,t)=y^*(x,t)$, we have
\begin{align}\label{3.2}
\hat{y}(x,t)=x.
\end{align}

Now we prove $\rho^{\tau}(x,t)$ defined in \eqref{1.3} weakly converges to $\rho_0(x)$ as $\tau\rightarrow 0$.
Firstly, by slow time scaling, 
\begin{align}
m^{\tau}(x,t):=m(x,\frac{t}{\tau})=
\begin{cases}
\int_{0{+}0}^{y_*^{\tau}(x,t){-}0}\,{\rm d}m_0(\eta), \quad{\rm if}\ \nu(x,t)=F^{\tau}(y_*^{\tau}(x,t);x,t),\\[1mm]
\int_{0{+}0}^{y_*^{\tau}(x,t){+}0}\,{\rm d}m_0(\eta), \quad{\rm if}\ \nu(x,t)=F^{\tau}(y_*^{\tau}(x,t){+}0;x,t);
\end{cases}
\end{align}
Then for fixed $t$ and almost all $x\in \mathbb{R}$, by \eqref{3.2}, we have
\begin{align}
\lim_{\tau\rightarrow0}m^{\tau}(x,t)=m_0(x),
\end{align}
which implies the weak convergence from $\rho^{\tau}(x,t)=m^{\tau}_x(x,t)$ to $\rho_0(x)$.

From \eqref{2.21b}, 
\begin{align}
\lim_{\tau\rightarrow0}\frac{1}{\tau}q^{\tau}(x,t):=\lim_{\tau\rightarrow0}\frac{1}{\tau}q(x,\frac{t}{\tau})
=\lim_{\tau\rightarrow0}\int_{0{+}0}^{y_*(x,\frac{t}{\tau})}u_0(\eta)\frac{1}{\tau}e^{-\frac{t}{\tau}}\,{\rm d}m_0(\eta)=0,
\end{align}
which, by \eqref{1.3b}, yields 
\begin{align}
\lim_{\tau\rightarrow0}u^{\tau}(x,t)&=\lim_{\tau\rightarrow0}\lim_{x_2,x_1\rightarrow x{\pm}0}\frac{1}{\tau}\frac{q^{\tau}(x_2,t)-q^{\tau}(x_1,t)}{m^{\tau}(x_2,t)-m^{\tau}(x_1,t)}\nonumber\\[1mm]
&=\lim_{x_2,x_1\rightarrow x{\pm}0}\lim_{\tau\rightarrow0}\frac{1}{\tau}\frac{q^{\tau}(x_2,t)-q^{\tau}(x_1,t)}{m^{\tau}(x_2,t)-m^{\tau}(x_1,t)}
=0.
\end{align}

\subsection{Vanishing--damping limit.}
Firstly, it follows from the formula of $\bar{F}(y;x,t)$ and $F(y;x,t)$ in \eqref{2.4} and \eqref{2.19} respectively that,
\begin{align}
\lim_{\tau\rightarrow\infty}F(y;x,t)=\bar{F}(y;x,t).
\end{align}
For $(x,t)$ satisfying $y_*(x,t)=y^*(x,t)$, by \eqref{2.21},
\begin{align}
x-t U\leq \liminf_{\tau\rightarrow\infty} y_*(x,t)\leq \limsup_{\tau\rightarrow\infty} y_*(x,t)\leq x+t U,
\end{align}
which implies that, there exists a convergent subsequence of $\{y_*(x,t)\}$ with respect to $\tau$ (still denoted as $y_*(x,t)$).
We denote the limit of the subsequence as
\begin{align}
\bar{Y}(x,t):=\lim_{\tau\rightarrow \infty}y_*(x,t)\in {\rm spt}\{\rho_0\}.
\end{align}
Similar to \eqref{3.1a}, we have 
\begin{align}
\bar{Y}(x,t)\in \bar{S}(x,t)\cap {\rm spt}\{\rho_0\}.
\end{align}
For fixed $t$ and almost all $x\in \mathbb{R}$, it follows from \eqref{2.10} that 
\begin{align}
\bar{Y}(x,t)=y_*(x,t),
\end{align}
which yields that
\begin{align}
\lim_{\tau\rightarrow\infty}m(x,t)=\bar{m}(x,t), \qquad \lim_{\tau\rightarrow\infty}q(x,t)=\bar{q}(x,t).
\end{align}
Thus the weak convergence from the measure $\rho(x,t)$ to $\bar{\rho}(x,t)$ is proved.

From \eqref{2.11} and \eqref{2.23},
\begin{align}
\lim_{\tau\rightarrow\infty}u(x,t)&=\lim_{\tau\rightarrow\infty}\lim_{x_2,x_1\rightarrow x{\pm}0}\frac{q(x_2,t)-q(x_1,t)}{m(x_2,t)-m(x_1,t)}\nonumber\\[1mm]
&=\lim_{x_2,x_1\rightarrow x{\pm}0}\lim_{\tau\rightarrow\infty}\frac{q(x_2,t)-q(x_1,t)}{m(x_2,t)-m(x_1,t)}=\bar{u}(x,t).
\end{align}

\smallskip
Up to now, Theorem 1.1 is proved.

\end{document}